\documentclass[twoside,10pt]{article}
\usepackage{cite,amsfonts,a4wide}
\usepackage{amsmath}
\usepackage{geometry}
\input amssym.def
\input amssym.tex

\newcommand{\BB}[1]{\mbox{${\Bbb #1}$}}

\newcommand {\R} {\mbox{$R$}}

\newcommand {\be}{\begin{equation}}
\newcommand {\e}{\end{equation}}
\newcommand {\bes}{\begin{displaymath}}
\newcommand {\es}{\end{displaymath}}
\newcommand {\nit}{\noindent}

\newcommand {\bit}{\bibitem}

\newcommand {\MF}{{\rm F}}

\newcommand {\MFO}{{\rm FO}}
\newcommand {\MFE}{{\rm FE}}
\newcommand {\MLW}{{\rm LW}}
\newcommand {\MUW}{{\rm UW}}
\newcommand {\MZC}{{\rm ZC}}
\newcommand {\MBR}{{\rm BR}}
\newcommand {\MWP}{{\rm WP}}
\newcommand {\MFT}{{\rm FT}}
\newcommand {\MWT}{{\rm WT}}
\newcommand {\MZR}{{\rm ZR}}
\newcommand{\Mfib}{{\rm fib}}
\newcommand{\Modfib}{{\rm odfib}}
\newcommand{\Mevfib}{{\rm evfib}}

\begin{document}


\begin{center}
{\bf \Large On the fibbinary numbers and the Wythoff array \\
\vspace{0.5cm}
A.J. Macfarlane$\;$}
\footnote{{\it e-mail}: 
a.j.macfarlane@damtp.cam.ac.uk}\\
\vskip .5cm
\begin{sl}
Centre for Mathematical Sciences, D.A.M.T.P.\\
Wilberforce Road, Cambridge CB3 0WA, UK
\vskip .5cm
\end{sl}
\end{center}

\begin{abstract}

This paper defines the set $\Mfib$ of fibbinary numbers and displays its
structure in the form of a table of a specialised type, and in array form.
It uses the Zeckendorf representation of
$n \in \BB{N}$ to define a bijection $\mathcal{Z}$
between $\BB{N}$ and $\Mfib$. It is
proved that the fibbinary array is the image under $\mathcal{Z}$ of the famous
Wythoff array. The fibbinary table provide useful pictorial insight into the
fractal defined by the Wythoff array. The Wythoff table, obtained as the image
under the inverse of $\mathcal{Z}$ of the fibbinary table, leads to another
simpler view of the fractal, and may be compared with the (1938) Steinhaus tree.

\end{abstract}

$\quad$  

$\quad$

\nit {\bf Keywords}: 
Fibonacci number, Zeckendorf representation, fibbinary number,
fibbinary array, Wythoff array, fractal sequence

$\quad$
 
\nit {\bf Mathematical Subject Classification}: 11A67, 11B83

\section{Introduction}

The purpose of this paper is to  define the set $\Mfib$ of fibbinary numbers,
$\underline{A003714}$ in \cite{oeis}, \cite{lsw},
present its integer elements

$\qquad\qquad$ (i) in table form as the fibbinary table, and

$\qquad\qquad$ (ii) in array form as the fibbinary array,
\vspace{0.5cm}
$\nit$ and to establish the following result.
\newline This states that the fibbinary array
is the image of the fascinating and widely studied Wythoff array,
$\underline{A035513}$ in \cite{oeis}, \cite{CS}, under the map
\begin{align} {\mathcal{Z}}: &\;\; \BB{N} \rightarrow \Mfib \nonumber \\
 & \;\; n \mapsto \Mfib(n). \label{aa0} \end{align} \nit
 This map, which depends upon the Zeckendorf representation \cite{as,zeck,lek}
 of the positive integers, is described in detail in
 Section three below.
 Tables like the fibbinary table play a valuable role below. They offer
 a clear view 
 of the structure not only of $\Mfib$ but also of the Wythoff fractal $W$
 \cite{ckFS, ckNS},
 which is defined by the Wythoff array. The Wythoff table is defined
 as the inverse image
 under $\mathcal{Z}$ of the fibbinary table. Of interest in its own right,
 it gives another view of the path from the Wythoff array
 to  its fractal $W$, and bears a close relationship 
 to the Steinhaus tree \cite{hs}.
 
  Section two is devoted to the Wythoff array, Table 1.
  The aim of Section two is to
 contrast the sheer weight of interesting material associated with
 it that has accumulated
 over the years with the transparent simplicity of its image under
 $\mathcal{Z}$, the fibbinary array.
 Section three treats the Zeckendorf representation \cite{as,zeck,lek} of the
 positive integers and uses it to define the fibbinary numbers,
 the bijection ${\mathcal{Z}}$, and the subset $\Modfib$ of $\Mfib$ of odd
 fibbinary numbers. 
 Section four introduces the fibbinary table, Table 2, and the
 fibbinary array, Table 3. It
 concludes with the statement and proof of the main result of this paper
 that the
 fibbinary array is the image under ${\mathcal{Z}}$ of the Wythoff array.
 Section five  discusses the Wythoff fractal, see Table 4, and introduces the
 Wythoff table, table 5.

\section{The Wythoff array}

The Wythoff array $\underline{A035513}$ in \cite{oeis}
is displayed, in extended form, in Table 1.
\vskip 0.5cm
\begin{center}
\begin{tabular}{|r|r||r|r|r|r|r|r|r|}
\hline
0 & 1 &    1 & 2  & 3  & 5   & 8   & 13  & 21    \\
1 & 3 &    4 & 7  & 11 & 18  & 29  & 47  & 76    \\
2 & 4 &    6 & 10 & 16 & 26  & 42  & 68  & 110   \\
3 & 6 &    9 & 15 & 24 & 39  & 63  & 102 & 165   \\
4 & 8 &   12 & 20 & 32 & 52  & 84  & 136 & 220   \\
5 & 9 &   14 & 23 & 37 & 60  & 97  & 157 & 254   \\
6 & 11 &  17 & 28 & 45 & 73  & 118 & 191 & 309   \\
7 & 12 &  19 & 31 & 50 & 81  & 131 & 212 & 343   \\
8 & 14 &  22 & 36 & 58 & 94  & 152 & 246 & 398   \\
9 & 16 &  25 & 41 & 66 & 107 & 173 & 280 & 453   \\
10 & 17 & 27 & 44 & 71 & 115 & 186 & 301 & 487   \\
11 & 19 & 30 & 49 & 79 & 128 & 207 & 335 & 542   \\
12 & 21 & 33 & 54 & 87 & 141 & 228 & 369 & 597   \\
\hline
\end{tabular}
\end{center}
\centerline{Table 1: The extended Wythoff array}
\vskip 0.5cm

The Wythoff array itself lies to the right of the double lines in Table 1.
A large amount of detailed information about it  is found under its entry in
\cite{oeis}. Also \cite{CS} contains a list of its seven most notable
properties, termed there, 
not unreasonably, as wonderful, properties 
which characterise it as an interspersion.
See, in this context, $\underline{A191426}$ in
\cite{oeis}, and \cite{ckID}. An additional property, similar in spirit to
the ones listed is described briefly in Section 2.1 below.
Of great interest too is the fact that the
Wythoff array defines a fractal sequence \cite{ckFS, ckNS}: see Section five.

The association of Wythoff's name with the Wythoff array has its origin in the
paper \cite{WAW} which describes a modification of the game of nim now referred
to as Wythoff's game \cite{HSMC}, \cite[p. 39]{WWRB}. A major role in the
win-lose strategy
of Wythoff's game is played by Wythoff pairs, $\MWP$,
\begin{align} \MWP = & \{ (x,y), \quad x=\MLW(n),\quad y=\MUW(n), \quad n=1,2,
\dots, \} \nonumber \\
= & \{ (1,2), (3,5), (4,7), (6,10), (8,13), (9,15), (11,18), (12,20),
\dots \} \label{A1} \end{align} \nit
Here the lower Wythoff sequence $\MLW$, $\underline{A000201}$ in \cite{oeis},
is
\begin{align} \MLW  = & \{ \MLW(n)=\lfloor n\phi \rfloor, \quad n=1,2,\dots \}
\nonumber \\
    = & \{ 1,3,4,6,8,9,11,12,14,16,17,19,21,22,24,25,27,29, \nonumber \\
    & 30,32,33,35,37,38,40,42,43,45,46,48,50,51,53,55, \dots \}
    \label{A2}, \end{align}\nit
and
the upper Wythoff sequence $\MUW$, $\underline{A001950}$ in \cite{oeis},
is
\begin{align} \MUW  = & \{ \MUW(n)= \lfloor n\phi^2 \rfloor,
\quad n=1,2,\dots \} \nonumber \\
    = & \{ 2,5,7,10,13,15,18,20,23,26,28,31,34,36,39,41, \nonumber \\
 & 44,47,49,52,54,57,60,62,65,68,70,73,75,78,81,\dots \}, \label{A3}
\end{align}\nit
where $\phi$ is the golden ratio, and
 $\lfloor x \rfloor$ is the floor function
\cite[p. 67]{gkp}, that is, $\lfloor x \rfloor$ is the greatest integer
less than or equal to $x, \; x \in \R$. Now various comments can be made.
First, the Wythoff array can be viewed as an array
of Wythoff pairs. Second, from (\ref{A1}), it can be seen that the integers
of $\MLW$ and $\MUW$  respectively appear in the columns with odd and even
column labels $k, \;\; k \geq 1$.  Third, the results
\be\label{A4}  \MLW \cup \MUW =\BB{N}, \quad \MLW \cap \MUW =\emptyset, \e\nit
agree with the property \cite{CS} that each $n \in \BB{N},
\;\; n \geq 1$, occurs exactly once in the Wythoff array. Also (\ref{A4})
implies that $\MLW$ and $\MUW$  provide an example of a pair of Beatty
sequences \cite{BeaPr,BeaS,Strut}.

\subsection{Fibonacci shapes and their covering property}

The Fibonacci shapes $\Gamma_n, \;\; n=1,2,  \dots$, are a subset of the
well-known set of Young shapes, with the definition:

The shape $\Gamma_n$ is a shape with successive columns of boxes of length
\be\label{A5} F_n,F_{n-1}, \dots, F_2,F_1, \quad n=1,2,\dots, \e\nit
and hence with $F_{n+2}-1$ boxes altogether.

The covering property of the Fibonacci shape $\Gamma_n$ is realised by placing
it over the top left-hand corner of, for example, the Wythoff array, thereby
covering a total of $F_{n+2}-1$ entries of the array. Then, for
$m=1,2,\dots$, those entries of the Wythoff array covered by $\Gamma_m$,
but not covered by $\Gamma_{m-1}$, are a set of $F_m$ integers $k$ such that
$F_{m+1} \leq k < F_{m+2}$, identified as Fibonacci subsets of $\BB{N}$
in Section 3.1.

\subsection{Construction and related facts}

There are two main ways of constructing the Wythoff array. First method: if the
columns of the array itself are labelled $k=1,2,\dots$, then the numbers
$n=0,1,2, \dots$ in the $k=-1$ column label the rows, and
the integers of the lower Wythoff sequence $\MLW$ of (\ref{A2})
provide the $k=0$ column.
Given the two columns to the left of the double line, the rest of the array is
obtained by using the Fibonacci addition rule. The $n=0,\;1$ rows contain the
Fibonacci and Lucas numbers. All the rows, and the columns with labels
$k=1,2, \dots, 6$, of the array displayed
in Table 1 have their own entries in \cite{oeis}.
An integer $n$ is said to be Fibonacci-odd, $\MFO$, or Fibonacci-even, $\MFE$, 
according as the LSB $a_2$ of its ZR $z(n)$ is odd or even. The sequence
$\MFO$, $\underline{A003622}$  in \cite{oeis}, is defined by
\begin{align} \MFO  = & \{ \MFO(n)=\lfloor (n+1)\phi^2 \rfloor-1,
\quad n=0,1,\dots \} \nonumber \\
= & \{1,4,6,9,12,14,17,19,22,25,27,30,33,35,38,40, \nonumber \\
& 43,46,48,51,53,56,59,61,64,67,69,72,74,77,80,
\dots \}.\label{A8} \end{align} \nit
The sequence $\MFE$, \underline{A022342}, in \cite{oeis} is defined, for
use in Section four, by
\begin{align} \MFE  = & \{ FE(n)=\lfloor (n+1)\phi \rfloor-1,
\quad n=0,1,\dots \} \nonumber \\
 = & \{ 0,2,3,5,7,8,10,11,13,15,16,18,20,21,23,24,26,28, \nonumber \\
 & 29,31,32,34,36,37,39,41,42,44,45,47,49,50,52,54,55, \dots \}
 \label{A7} \end{align}\nit

The integers of $\MFO$ provide the first or $k=1$ column of the Wythoff array.
Further, if $C_k(n)$ is the element of the Wythoff array in
row $n$ and column $k$,
then \cite{CS}
\be\label{A11} C_k(n)= \lfloor (n+1) \phi \rfloor F_{k+1}+nF_k,\e\nit
so that for the first or $k=1$ column
\be\label{A12} C_1(n)=\lfloor (n+1) \phi \rfloor +n 
=\lfloor (n+1) \phi^2 \rfloor-1 = \MFO(n). \e\nit

A further property of the columns of the Wythoff array is given in
Section 3.1 after the Zeckendorf representation of the
positive intgers has been defined.

Second method: given the first row of the array itself, then \cite{ewwWA} build
subsequent rows by iteratively adding
\be\label{A9} \{ F_{3+s}, \; F_{4+s}, \; F_{5+s}, \; \dots \}, \e\nit
where $s=0$ or $1$ is the smallest offset producing an initial term that
has not occurred in an earlier row. Rows with $s=0, 1$ respectively are
numbered \cite{ewwWA} by the integers of the upper and lower Wythoff
sequences: see (\ref{A3}) and (\ref{A2}).

Relationships between the sequences $\MFO, \MFE, \MLW, \MUW$ include
\begin{align} \MFE(n-1)  = & \MLW(n)-1, \quad n=1,2,\dots \nonumber \\
\MFO(n-1)  = & \MUW(n)-1, \quad n=1,2,\dots \nonumber \\
\MUW(n)  =  & \MLW(n) + n, \quad n=1,2,\dots. \label{A10} \end{align} \nit

\section{The Zeckendorf representation and the fibbinary numbers}
\subsection{The Zeckendorf representation}

Let $\MF$ denote the set of Fibonacci numbers
\begin{align} \MF & =\{ F_n, \quad n \geq 2 \} \nonumber \\
& = \{1,2,3,5,8,13,21,34,55, \dots \}. \label{aa1} \end{align} \nit

Zeckendorf's theorem \cite{zeck,lek,as} states that every
positive integer can be expressed uniquely as a sum of
non-consecutive Fibonacci numbers $F_n \in \MF$.
This leads to the Zeckendorf representation -- notation $z(n)$ and
abbreviation ZR -- of the integers $n \in \BB{N}, \; n \geq 1$.
This 
can be written uniquely in the form
\be \label{aa2} z(n)= \sum_{i=2}^r a_i F_i, \e\nit
where
\begin{align}
& a_i \in \{ 0,1 \},   \quad a_r=1, \nonumber \\
  & a_i\;a_{i+1}=0,  
\quad 2 \leq i <r, \label{A13}  \end{align} \nit
and $r$ is the largest integer such that $F_r \leq n$. The entry $a_2$ is
referred to as the least significant bit, LSB, and (\ref{A13}) as the
Zeckendorf condition, $\MZC$.

For $z(n)$, given  by (\ref{aa2}), write also
\be\label{aa3} z(n)=(a_r,a_{r-1}, \dots, a_2)_F = (n)_F, \e\nit
where, within $(n)_F$, $n$ denotes the
string
\be\label{aa4} a_r a_{r-1} \dots a_2, \e\nit
of ones and zeros obtained from (\ref{aa2}). Define also the length $l(n)$
of the string (\ref{aa4}), that is, the total number of ones and zeros.
A formula for $l(n)$ , due to R. Stephan, is found under the entry
\underline{A072649} in \cite{oeis}.

Note some examples
\begin{align} z(F_n)= (F_n)_F= (10^{n-2})_F, \quad & l(F_n)=n-1 \nonumber \\
z(F_{2r}-1)= ( (10)^{r-1})_F,\quad & l(F_{2r}-1)=2r-2 \nonumber \\
z(F_{2r+1}-1)= ( (10)^{r-1}1)_F,\quad & l(F_{2r+1}-1)=2r-1 \label{aa5}.
\end{align}\nit

Also the $k$-th column of the Wythoff array contains only those integers
$n \in \BB{N}$
for which the string $(n)_F$ defined by (\ref{aa4}) terminates $10^{k-1}$.
For example, the integers $8,29,42 $ of column $k=5$ have Zeckendorf
representation $z(n)=(n)_F$ with string $n$ given by
$10^4,\; 1010^4,\; 10010^4$.

If the integers $n \in \BB{N}$ are listed in ZR form, it is seen that they
separate into Fibonacci subsets ${\mathcal{F}}_k$ of cardinality $F_k$ such
that $l(n)=k$ for all $n \in {\mathcal{F}}_k$. Thus, for $k=6$,
\be\label{aa6} {\mathcal{F}}_6= \{ 13,14, \dots, 20 \},\;  F_6\;=8, \e\nit
\be\label{aa7} z(13)=(10^5)_F, \;z(14)=(10^41)_F,\; \dots, \;
z(20)=((10)^3)_F,\e\nit
and $13=F_7, 20=F_8-1$.

The sequence $\MFO$ of Fibonacci-odd integers is defined by (\ref{A8}) in
section 2.2.
To prove a result
needed below, namely
\be\label{aa21} z(FO(n))= (n01)_F, \e\nit
where $n$ is the string given as  $(n)_F$ by (\ref{aa3}),
set out from (\ref{A8}). Thus
\begin{align} \MFO(n)= & \lfloor(n+1)\phi^2 \rfloor-1=
\lfloor(n+1)\phi \rfloor +n \nonumber \\
= & (n0)_F+n+1 \nonumber \\
= & (n00)_F+1 \nonumber \\
= & (n01)_F. \label{aa22} \end{align}\nit
The second step uses a result proved in \cite{Rebie}
\be\label{aa23}
(n0)_F=\lfloor (n+1) \phi \rfloor-1.\e\nit The third 
step uses the Fibonacci addition rule. 

\subsection{The fibbinary numbers}

The definition of $z(n)$ of (\ref{aa2}) and (\ref{aa3}) motivates
introduction of the set
\be\label{aa25} \Mfib = \{\Mfib(n), \quad n=1,2, \dots \}, \e\nit
$\underline{A003714}$ in \cite{oeis},\cite{lsw},
of fibbinary numbers, wherein the positive integer $\Mfib(n)$ is defined by
\be\label{bb7} \Mfib(n)=\sum_{i=2}^r a_i 2^{i-2}, \e\nit
and use is made of the same set of
coefficients $a_i, \quad 2 \leq i < r$,
as are used in (\ref{aa2}). Eq. (\ref{bb7}) is the binary representation,
$\MBR$, of $\Mfib(n)$. Also, using (\ref{aa5}), it follows that for
$F_i \in \MF$
\be\label{bb9} \Mfib(F_i)=2^{i-2}, \quad i \geq 1.\e\nit
so that the $\Mfib(F_i)$ provide a basis for the fibbinary numbers.

Just as for the integers $n \in \BB{N}$ in ZR form,
the $\Mfib(n)$
separate into Fibonacci subsets $\mathcal{F}_k$ of cardinality $F_k$ with
$l(\Mfib(n))=l(n)$. Accordingly the listing of the elements of $\Mfib$
is presented as
\begin{align}
\Mfib  = & \{\Mfib(n), \quad n=1,2,\dots \} \nonumber \\
 =  & \{ 1;2;4,5;8,9,10;16,17,18,20,21=\Mfib(12);32,33,34,36,37,40,41,
 42=\Mfib(20); \nonumber \\
 & 64,65,66,68,69,72,73,74,80,81,82,84,85=\Mfib(33); \nonumber \\
 & 128, \dots , 170=\Mfib(54); 256,\dots, 341=\Mfib(88); \dots \}, \label{b11}
\end{align}\nit
using semi-colons to indicate separation of the
entries into Fibonacci subsets. The largest element of the
$\mathcal{F}_k, \;\; k=1,2,\dots $,
subset of $\Mfib$ is $\Mfib(F_{k+2}-1)$, as noted explicitly in (\ref{b11})
for $k=5, \dots, 9$.

\nit Writing (\ref{bb7}) as
\be\label{bb8} \Mfib(n)=(a_r,a_{r-1} \dots a_2)_2 =(n)_2, \e\nit
in analogy with (\ref{aa3}), see the same sequence (\ref{aa4}) used for $n$
in $(n)_F$ and $(n)_2$ but with distinct meanings. Respect of the Zeckendorf
condition is a strict requirement in both uses.

It is clear that the map (\ref{aa0}) from  $\BB{N}$ to $\Mfib$
defines a bijection
which may be expressed now more explicitly as
\begin{align} {\mathcal{Z}}: & \;\; \BB{N} \rightarrow \Mfib \nonumber \\
 & \;\; n =(n)_F \mapsto \Mfib(n)=(n)_2 \label{bb12}. \end{align} \nit
For example, for $n=12$,
\begin{align} z(12)= & (10101)_F=F_6+F_4+F_2=8+3+1=12 \nonumber \\
\Mfib(12)= & (10101)_2= 16+4+1=21. \label{bb13} \end{align}\nit
 It may be checked that the places the entries $12$ and $21$ occupy in
 Tables 1 and 3 agree.

\subsection{The key property of $\Mfib$ and the set $\Modfib$}

Define next the important subset $\Modfib$ of $\Mfib$ of odd fibbinary numbers,
$\underline{A022341}$ in \cite{oeis}, \cite{lsw},
\begin{align} \Modfib  = &  \{ \Modfib(n), \quad n=0,1,\dots \}
\nonumber \\
 = & \{1,5;9;17,21;33,37,41;65,69,73,81,85; \nonumber \\
 & 129,133,137,145,149,161,165,169; \nonumber \\
 & 257, \dots ,341; 513, \dots \} \label{bb14}
\end{align}\nit
where $\Modfib(0)=1$ requires the use of $\Mfib(0)$.
The even fibbinary numbers of $\Mevfib$, $\underline{A022342}$ in
\cite{oeis}, are not used here.

The set $\Mfib$ has the key property that
if any integer $j$ occurs in $\Mfib$, then so also does the even integer $2j$
and the odd integer $4j+1$,
but no integer of the form $4j-1$ does.

It follows from the key property of $\Mfib$ that
\begin{align} \Modfib(n) & = 4\Mfib(n)+1 =(n01)_2, \; n=0,1, \dots,
\label{bb15} \\
\Mevfib (n) &= 2\Mfib(n)=(n0)_2, \; n=1,2, \dots.\nonumber \end{align} \nit

It is natural to expect that $\MFO$ and $\Modfib$ are closely related.
Eq. (\ref{aa21}) states that $\MFO(n)$ has ZR
\be \label{bb16} z(\MFO(n))= (n01)_F \e\nit
so that its image under $\mathcal{Z}$ is
\be\label{bb20} \Mfib(\MFO(n))= (n01)_2=\Modfib(n). \e\nit
Thus there is a bijection between $\MFO$ and $\Modfib$ inherited from
(\ref{bb12}).

It may also be noted that the odd elements of the $\mathcal{F}_k$ subsets of
$\Mfib$  are equal to the elements of the $\mathcal{F}_{k-2}$ subsets of
$\Modfib$ for $k=3,4,\dots$.

\section{The fibbinary table and the fibbinary array}

Setting out from the lowest element $\Mfib(1)=1$ of $\Mfib$, and
using the key property of $\Mfib$, it is easy to generate all the integer
entries of $\Mfib$. The way of displaying them preferred here involves drawing
up a table like Table 2.

\vspace{0.5cm}
\begin{center}
\begin{tabular}{|r|r|r|r|r|r|r|r|r|r|r|r|r|}
\hline
1   & & & & &  & & & &    & & &  \\
2   & & & & &  & & & &    & & &    \\
4   & & & & &  & & & 5    & & & & \\
8   & & & & &  9 & & & 10 & & & &   \\
16  & & & 17 & & 18 & & & 20 & & & 21 &   \\
32  & & 33 & 34 & & 36 & & 37 & 40 & & 41 & 42 &  \\
\hline
\rule[-2mm]{0mm}{6mm} 0 & 1 & 2 & 3 & 4 & 5 & 6 & 7 &
8 & 9 & 10 & 11 & 12   \\
\hline
\end{tabular}
\end{center}
\vskip 0.5cm
\centerline{Table 2: The fibbinary table, $\MFT$}
\vskip 0.5cm
Various features of such tables are seen clearly in Table 2. First, at the head
of its columns, find  the integers
\be\label{ss1} 1,5,9,17,21,33,37,41, \dots, \e\nit
of the sequence $\Modfib$ of (\ref{bb14}). Second, below each one of these,
find all its even multiples one after another.
Third, the elements of the
Fibonacci subsets $\mathcal{F}_k, \;\; k=1,\dots,6$, of $\Mfib$ are aligned
along rows $k=1,\dots,6$.

Tables like Table 2 are drawn up by the following procedure. Suppose a table
with $p$ rows $p=1,2, \dots$ is called for. First, use the labels
\be\label{ss2}  x=0,1, \dots, F_{p+1}-1, \e\nit
to provide  the required number of columns. Second, to create the
bottom, or $p$-th, row fill the
integers of the Fibonacci subset $\mathcal{F}_p$ in order into the columns
\be\label{ss3} x=0,2,3,5,7,8,10,11,13, \dots, x_{max} , \e\nit
where $x_{max}$ is to be determined. Third, find all the even entries of
$\mathcal{F}_p$, and, directly above them, place entries obtained by dividing
by two, in each case repeating the process until an odd number is reached.
The numbers in (\ref{ss3}) belong to the sequence $\MFE$ of (\ref{A7}), and
can therefore be written as
\be\label{ss4}  \MFE(n), \quad n=0,1, \dots F_p-1, \e\nit
since $\mathcal{F}_p$ has $F_p$ elements. Hence
\be\label{ss5} x_{max}= \MFE(F_p-1). \e\nit
The empty columns of the $p$-th row of the $p$ row table have labels drawn
from the sequence $\MFO$, (\ref{A8}).

For $p=6$ it is straightfroward to confirm that the procedure described produces
Table 2.

The features such as those noted in Table 2, give, first, the following
description of
the structure of the set $\Mfib$
\begin{align} \Mfib= & \bigcup_{s=0}^\infty \alpha(s),  \quad
\alpha(s)\bigcap \alpha(s^{\prime})=\emptyset, \quad s \neq s^{\prime},
\nonumber \\
\alpha(s)= & \Modfib(s) \{2^r, \, r=0,1,\dots, \} \label{f1} \end{align}\nit
and, second, motivates the display of the integers of $\Mfib$ in array form. 
This has the integers $\Modfib(n)$ of
$\Modfib$, taken in the order $n=0,1,\dots$, in its 
first column and their multiples by powers of $2$ in the subsequent columns.
Thus the fibbinary array takes the form shown in Table 3.
\begin{center}
\begin{tabular}{|r|r|r|r|r|r|r|}
\hline
1 & 2 & 4 & 8 & 16 & 32 & 64    \\
5 & 10 & 20 & 40 & 80 & 160 & 320  \\
9 & 18 & 36 & 72 & 144 & 288 & 576  \\
17 & 34 & 68 & 136 & 272 & 544 & 1088   \\
21 & 42 & 84 & 168 & 336 & 672 & 1344  \\
33 & 66 & 132 & 264 & 528 & 1056 & 2112\\
37 & 74 & 148 & 296 & 592 & 1194 & 2388 \\
41 & 82 & 164 & 328 & 656 & 1312 & 2624 \\
65 & 130 & 260 & 520 & 1040 & 2080 & 4160 \\
69 & 138 & 276 & 552 & 1104 & 2208 & 4416 \\
73 & 146 & 292 & 584 & 1168 & 2336 & 4672 \\
81 & 162 & 324 & 648 & 1296 & 2592 & 5184 \\
85 & 170 & 340 & 680 & 1360 & 2720 & 5440 \\
\hline
\end{tabular}
\end{center}
\vskip 0.5cm
\centerline{Table 3: The fibbinary array}
\vskip 0.5cm

The statement that the fibbinary array is the image under the bijection
${\mathcal{Z}}$ of (\ref{aa0}) and (\ref{bb12}) of the Wythoff array is the 
main result of this paper.
It is known already that the image  under ${\mathcal{Z}}$ of the entry
$\MFO(n)$ in $n$-th place in
the first column of the
Wythoff array is $\Modfib(n)$ in $n$-th place in
the first column of the fibbinary array. To complete proof of the statement
made requires proof that
\be\label{f2} \Mfib(C_k(n))= 2^{k-1} \Modfib(n), \; k=1,2,\dots, \;
n=0,1,\dots \; .\e\nit

\subsection{Proof of (\ref{f2})}

Set out from the result from \cite{CS} given above as (\ref{A11}), and apply
once again (\ref{aa23}) to get
\begin{align} C_k(n)= & \lfloor (n+1)\phi \rfloor F_{k+1}+nF_k \nonumber \\
= & ((n0)_F+1) F_{k+1}+nF_k. \label{f11} \end{align} \nit
Now use the result, an easy generalisation of the Fibonacci addition rule, 
\be\label{f12} F_{k+1}(n0)_F+F_k(n)_F=(n0^{k+1})_F, \e\nit
which may be proved by induction on $k$. This leads, using (\ref{aa5}), to 
\be\label{f13} C_k(n)= (n0^{k+1})_F+F_{k+1}= (n010^{k-1})_F. \e\nit
Hence, using (\ref{bb20}),
\begin{align} \Mfib(C_k(n))= & (n010^{k-1})_2 \nonumber \\
= & 2^{k-1}(n01)_2 =2^{k-1} \Modfib(n). \label{f14} \end{align} \nit
This completes the proof that the image under ${\mathcal{Z}}$ of the
$(n,k)$ element of the Wythoff array is the $(n,k)$ element of the
fibbinary array.

\section{The Wythoff fractal }

The section begins with a brief review of well-known material on the Wythoff
fractal, and goes on to consider  it from the viewpoint of this paper.
In this section and only here, the rows of
the Wythoff array are numbered
$n=1,2,\dots$, and the same applies to the fibbinary array.

A fractal sequence \cite{ckFS, ckNS}
is a sequence which contains itself as a proper subsequence
infinitely many times. The fractal sequence $W$ of the Wythoff array
\cite{ckFS, ckNS} has the
following definition:
\be\label{dd1} W= \{ w(n), \; n=1,2,\dots\} \e\nit
where $w(n)$ is the number of the row of the Wythoff array to which the
integer $n$ belongs.
This yields the sequence, $\underline{A003603}$ in \cite{oeis},
\be\label{dd2} W=\{ 1;1;1,2;1,3,2;1,4,3,2,5;1,6,4,3,7,2,8,5;1,\dots \}
\e\nit

The fractal nature of the fractal can be seen in a statement
which can be applied repeatedly:
when the first occurrence of each integer of the sequence is excised
what remains is identical to the original sequence.

The fibbinary array yields the same fractal in exactly the same way.

The fibbinary table, Table 2, can be used to get a tabular view of $W$.
If each integer entry of Table 2 is replaced by the number of the row of the
fibbinary array in which it is found, then the fractal table, Table 4, is 
obtained.

\vspace{0.5cm}
\begin{center}
\begin{tabular}{|r|r|r|r|r|r|r|r|r|r|r|r|r|}
\hline
1   & & & & &  & & & &    & & &  \\
1   & & & & &  & & & &    & & &    \\
1           & & & & &  & & & 2    & & & & \\
1     &   & & & & 3 & & &    2 & & & &   \\
1     & & & 4 & & 3 & & &    2 & & &   5 &   \\
1  & & 6 &  4 & & 3 & & 7 &  2 & & 8 & 5 &  \\
\hline
\rule[-2mm]{0mm}{6mm} 0 & 1 & 2 & 3 & 4 & 5 & 6 & 7 &
8 & 9 & 10 & 11 & 12   \\
\hline
\end{tabular}
\end{center}
\vskip 0.5cm
\centerline{Table 4: The fractal table}
\vskip 0.5cm

The
Fibonacci subsets $\mathcal{F}_k, \; k=1,\dots,6$, of $W$ occupy the rows of
Table 4. The numbers of the subset $\mathcal{F}_7$ of $W$, in correct fractal
order, are found from Table 4 simply by filling the integers $9=F_6+1, \dots,
13=F_7$ into the spaces in row six of the table. This yields
\be\label{dd22} 1,9,6,4,10,3,11,7,2,12,8,5,13. \e\nit
It may be checked, by drawing up a seven row table, that the numbers
(\ref{dd22}) appear in its seventh row.  Now, by filling the numbers
$14=F_7+1, \dots, 21=F_8$ in order into the spaces in that row, the numbers of
the $\mathcal{F}_8$ subset of $W$ can be read off.

\subsection{The Wythoff table, $\MWT$}

The Wythoff table, $\MWT$, is defined as the inverse image under $\mathcal{Z}$
of the fibbinary table of Table 2. $\MWT$ is presented in Table 5.
\vspace{0.5cm}
\begin{center}
\begin{tabular}{|r|r|r|r|r|r|r|r|r|r|r|r|r|}
\hline
1   & & & & &  & & & &    & & &  \\
2   & & & & &  & & & &    & & &    \\
3   & & & & &  & & & 4    & & & & \\
5  & & & & &  6 & & & 7 & & & &   \\
8 & & & 9 & & 10 & & & 11 & & & 12 &   \\
13  & & 14 & 15 & & 16 & & 17 & 18 & & 19 & 20 &  \\
\hline
\rule[-2mm]{0mm}{6mm} 0 & 1 & 2 & 3 & 4 & 5 & 6 & 7 &
8 & 9 & 10 & 11 & 12   \\
\hline
\end{tabular}
\end{center}
\vskip 0.5cm
\centerline{Table 5: The Wythoff table}
\vskip 0.5cm

As is to be expected, Table 4 may be obtained from Table 5, just as
it was from
Table 2. Nevertheless Table 5 has several interesting features. First, at the
head of its columns are found the integers of the sequence $\MFO$ of (\ref{A8}).
In view of (\ref{bb20}), this agrees with the fact that the integers of
$\Modfib$
are found at the head of the columns of Table 2. Second, Fibonacci addition
applies down the columns of Table 5. But there remains the question: what,
for $\MWT$, plays the role for $\MFT$ of the key property of $\Mfib$. To
answer it, use the $\MZR$ of $z(n)$ of (\ref{aa2}) of the
integer $n \in \BB{N}$.
Given an integer $j \in \BB{N}$ in  row $m$, with $z(j)=(j)_F$,
the even integers
\be\label{dd23} (j0)_F, (j00)_F, \dots, \e\nit
are found in the same column in rows $m+1$ and $m+2$. This agrees with the
fact that Fibonacci addition 
\be\label{dd4} (j00)_F=(j0)_F+(j)_F. \e\nit
applies down the columns of Table 5. Also the odd
integer with $\MZR, \; (j01)_F$ is found also in row $m+2$ , in the next
non-empty column to the right of $(j00)_F$.
For example, for $j=7$,
\begin{align} j= & 7=(1010)_F, \; (j0)_F=(10100)_F=11, \; (j00)_F=18, 
\nonumber \\
(j01)_F= & (101001)_F=11+7+1=19. \label{dd5} \end{align}

The information in the Wythoff Table could very well have been
displayed in the form of a tree. This would take the same form
as the Steinhaus tree \cite{hs}, which has
appeared previously, long ago in fact, in 1938, but not in the context of $W$.
See also \cite{TWvR}.

A final comment: whereas the fibbinary array is seen to be a much simpler
entity than the
Wythoff array, the Wythoff table appears to be simpler than the fibbinary table.

\end{document}